\documentstyle[11pt]{article}
\begin{document}

\noindent{\bf CENTRAL ELEMENTS OF THE ALGEBRAS $U'_q({\rm so}_m)$ AND
$U_q({\rm iso}_m)$}
\bigskip

\centerline{\sl M. Havl\'\i\v cek, S. Po\v sta}

\centerline{Department of Mathematics and Dopler Institute, FNSPE, CTU}

\centerline{Trojanova 13, CZ-120 00 Prague 2, Czech Republic}
\medskip

\centerline{\sl A. Klimyk}
\centerline{Institute for Theoretical Physics,
Kiev 252143, Ukraine}
\bigskip

\begin{abstract}
The aim of this paper is to give a set of central elements of the
algebras $U'_q({\rm so}_m)$ and $U_q({\rm iso}_m)$
when $q$ is a root of unity. They are surprisingly arise
from a single polynomial Casimir element of the algebra
$U'_q({\rm so}_3)$. It is conjectured that the Casimir elements
of these algebras under any values of $q$ (not only for $q$ a root of unity)
and the central elements for $q$ a root of unity derived in this paper
generate the centers of $U'_q({\rm so}_m)$ and $U_q({\rm iso}_m)$ when $q$
is a root of unity.

\end{abstract}


\newfam\eusfam
\font\tenseuler=msbm10
\textfont\eusfam=\tenseuler
\def\dvoj#1{{\fam\eusfam #1}}
\def\C{\dvoj C}
\def\N{\dvoj N}
\def\R{\dvoj R}
\newfam\eusfamdva
\font\tensdvaeuler=msam10
\font\sevsdvaeuler=msam7
\font\fivsdvaeuler=msam5
\textfont\eusfamdva=\tensdvaeuler
\scriptfont\eusfamdva=\sevsdvaeuler
\scriptscriptfont\eusfamdva=\fivsdvaeuler
\def\sixt#1{\ifcase#1 0\or1\or2\or3\or
4\or5\or6\or7\or8\or9\or A\or
B\or C\or D\or E\or F\fi}
\newbox\strutboxdva
\setbox\strutboxdva=\hbox{\vrule height2.5pt depth3.5pt
width 0pt}
\def\strutdva{\relax\ifmmode\copy\strutboxdva
\else\unhcopy\strutboxdva\fi}
\def\equivs{\equiv\discretionary{}{\hbox{$\equiv$}}{}}
\def\komb#1#2{\bigl(
\mskip-0.7\thinmuskip{{#1\atop #2}}
\mskip-0.7\thinmuskip\bigr)}
\def\krat{\!\cdot\!}
\def\kombd#1#2{\Bigl(
\mskip-0.7\thinmuskip{{#1\atop #2}}
\mskip-0.7\thinmuskip\Bigr)}
\mathchardef\leq="3\sixt\eusfamdva35
\mathchardef\geq="3\sixt\eusfamdva3D
\def\mod{\,{\rm mod}\,}

{\bf 1.} The algebra $U'_q({\rm so}_m)$ is a nonstandard $q$-deformation
of the universal enveloping algebra $U({\rm so}_m)$ of the Lie algebra
${\rm so}_m$.
It was defined in [1] as the associative algebra (with a unit)
generated by the elements $I_{21}$, $I_{32},\cdots ,I_{m,m-1}$
satisfying the defining relations
$$I_{i+1,i}I^2_{i,i-1}-(q+q^{-1})I_{i,i-1}
I_{i+1,i}I_{i,i-1}
+I^2_{i,i-1}I_{i+1,i}=-I_{i+1,i},\eqno (1)$$
$$I^2_{i+1,i}I_{i,i-1}-(q+q^{-1})I_{i+1,i}I_{i,i-1}
I_{i+1,i}+I_{i,i-1}I^2_{i+1,i}=-I_{i,i-1},\eqno(2)$$
$$[I_{i,i-1},I_{j,j-1}]=0\ \ \ {\rm for}\ \ \
\vert i-j\vert >1, \eqno (3) $$
where [.,.] denotes the usual commutator.
In the limit $q\to 1$
formulas (1)--(3) give the relations defining the universal enveloping
algebra $U({\rm so}_m)$. Note also that the relations
(1) and (2) differ from the $q$-deformed Serre relations
in the approach of Drinfeld and Jimbo to quantum orthogonal
algebras (see, for example, [2])
by presence of nonzero right hand sides in (1) and (2).

For the algebra $U'_q({\rm so}_3)$ the relations (1)--(3) are reduced to the
following two relations:
$$
I_{32}I^2_{21}-(q+q^{-1})I_{21}
I_{32}I_{21}+I^2_{21}I_{32}=-I_{32},\eqno (4)$$
$$I^2_{32}I_{21}-(q+q^{-1})I_{32}I_{21}
I_{32}+I_{21}I^2_{32}=-I_{21}. \eqno(5)
$$
Denoting $I_{21}$ and $I_{32}$ by $I_1$ and $I_2$, respectively, and
introducing the element $I_3:=
q^{1/2} I_1 I_2-q^{-1/2} I_2 I_1$ we find that relations (4) and (5) are
equivalent to three relations
$$q^{1/2} I_1 I_2-q^{-1/2} I_2 I_1=I_3,\eqno(6)$$
$$q^{1/2} I_2 I_3-q^{-1/2} I_3 I_2=I_1,\eqno(7)$$
$$q^{1/2} I_3 I_1-q^{-1/2} I_1 I_3=I_2.\eqno(8)$$

The Inonu--Wigner contruction applied to the algebra $U'_q({\rm so}_m)$
leads to the algebra $U_q({\rm iso}_{m-1})$ which was defined in [3].
The algebra $U_q({\rm iso}_{m})$ is the associative algebra (with a unit)
generated by the elements $I_{21}$, $I_{32},\cdots ,I_{m,m-1}$, $T_m$
such that the elements $I_{21}$, $I_{32},\cdots ,I_{m,m-1}$ satisfy
the defining relations of the algebra $U'_q({\rm so}_m)$ and the
additional defining relations
$$
I^2_{m,m-1}T_m-(q+q^{-1}) I_{m,m-1}T_mI_{m,m-1}+T_mI^2_{m,m-1}=-T_m,
$$
$$
I_{m,m-1}T^2_m-(q+q^{-1}) T_mI_{m,m-1}T_m+T^2_mI_{m,m-1}=-I_{m,m-1},
$$
$$
[I_{k,k-1},T_m]:=I_{k,k-1}T_m-T_mI_{k,k-1}=0\ \ \ \ {\rm if}\ \ \ \
k<m.
$$
If $q=1$, then these relations define the universal enveloping algebra
$U({\rm iso}_m)$ of the Lie algebra ${\rm iso}_m$ of the Lie group
$ISO(m)$.

The algebra $U({\rm iso}_2)$ is generated by two elements $I_{21}$ and
$T_2$ satisfying the relations
$$
I^2_{21}T_2-(q+q^{-1})I_{21}T_2I_{21}+T_2I^2_{21}=-T_2, \eqno (9)$$
$$
I_{21}T^2_2-(q+q^{-1})T_2I_{21}T_2+T^2_2I_{21}=-I_{21}. \eqno (10)
$$
Denoting $I_{21}$ by $I$ and introducing the element
$T_1:=[I,T_2]_q\equiv q^{1/2}IT_2-q^{-1/2}T_2I$, we find that the
relations (9) and (10) are equivalent to the relations
$$
[I,T_2]_q=T_1,\ \ \ [T_1,I]_q=T_2,\ \ \ [T_2,T_1]_q=0. \eqno (11)
$$

{\bf 2.} In $U'_q({\rm so}_m)$
we can determine [4] elements analogous to the basis
matrices $I_{ij}$, $i>j$, (defined, for example, in [5]) of ${\rm so}_m$.
In order to give them we use the
notation $I_{k,k-1}\equiv I^+_{k,k-1}\equiv I^-_{k,k-1}$.
Then for $k>l+1$ we define recursively
$$
I^{\pm}_{kl}:= [I_{l+1,l},I_{k,l+1}]_{q^{\pm 1}}\equiv
q^{\pm 1/2}I_{l+1,l}I_{k,l+1}-
q^{-\pm 1/2}I_{k,l+1}I_{l+1,l}. \eqno (12)
$$
The elements $I^+_{kl}$, $k>l$, satisfy the commutation relations
$$
[I^+_{ln},I^+_{kl}]_q=I^+_{kn},\ \
[I^+_{kl},I^+_{kn}]_q=I^+_{ln},\ \
[I^+_{kn},I^+_{ln}]_q=I^+_{kl} \ \ \
{\rm for}\ \ \  k>l>n, \eqno (13) $$
$$
[I^+_{kl},I^+_{nr}]=0\ \ \ \ {\rm for}\ \ \
k>l>n>r\ \ {\rm and}\ \ k>n>r>l, \eqno (14)
$$
$$
[I^+_{kl},I^+_{nr}]_q=(q-q^{-1})
(I^+_{lr}I^+_{kn}-I^+_{kr}I^+_{nl}) \ \ \ {\rm for}\ \ \
k>n>l>r. \eqno (15)
$$
For $I^-_{kl}$, $k>l$, the commutation relations are obtained by replacing
$I^+_{kl}$ by $I^-_{kl}$ and $q$ by $q^{-1}$.

Using the diamond lemma (see, for example, Chapter 4 in [2]),
N. Iorgov proved the Poincar\'e--Birkhoff--Witt theorem for the
algebra $U'_q({\rm so}_m)$ (proof of it will be published):
\medskip

\noindent
{\bf Theorem 1.} {\it The elements ${I_{21}^+}^{n_{21}}{I_{32}^+}^{n_{32}}
{I_{31}^+}^{n_{31}}\cdots {I_{m1}^+}^{n_{m1}}$, $n_{ij}=0,1,2,\cdots $,
form a basis of the algebra $U'_q({\rm so}_m)$.}
\medskip

This theorem is true if the elements $I^+_{ij}$ are replaced by the
corresponding elements $I^-_{ij}$.

Using the generating elements $I_{21}$, $I_{32},\cdots ,I_{m,m-1}$
of the algebra $U_q({\rm iso}_m)$ we define by formula (12) the
elements $I^{\pm}_{ij}$, $i>j$, in this algebra. Besides,
in $U_q({\rm iso}_m)$ we also define recursively the elements
$$
T^{\pm}_k:=[I_{k+1,k},T^{\pm}_{k+1}]_{q^\pm},\ \ \ k=1,2,\cdots ,m-1 .
$$
It is shown in [6] that
the elements $I^+_{ij}$, $i>j$, and $T^+_k$, $1\le k\le m$,
satisfy the commutation relations (13)--(15) and the relations
$$
[I^+_{ln},T^+_l]_q=T^+_n, \ \  [T^+_n, I^+_{ln}]_q=T^+_l \ \ \
{\rm for}\ \ \ l>n, $$
$$
[T^+_l,I^+_{np}]=0 \ \ \ \ {\rm for}\ \  l>n>p\ \ {\rm or}\ \ n>p>l,
$$
$$
[T^+_l,I^+_{np}]=(q-q^{-1})(T^+_nI^+_{lp}-T^+_pI^+_{nl})\ \ \
{\rm for} \ \ \ n>l>p,
$$
$$
[T^+_l,T^+_n]_q=0\ \ \ {\rm for}\ \ \ n<l.
$$
For $U_q({\rm iso}_m)$ the Poincar\'e--Birkhoff--Witt theorem
is formulated as
\medskip

\noindent
{\bf Theorem 2.} {\it The elements ${I_{21}^+}^{n_{21}}{I_{32}^+}^{n_{32}}
{I_{31}^+}^{n_{31}}\cdots {I_{m1}^+}^{n_{m1}}{T_1^+}^{n_1}{T_2^+}^{n_2}
\cdots {T_m^+}^{n_m}$ with $n_{ij},n_k=0,1,2,\cdots $,
form a basis of the algebra $U_q({\rm iso}_m)$.}
\medskip

{\bf 3.} It is easy to check that
for any value of $q$ the algebra $U'_q({\rm so}_3)$ has the
Casimir element
$$
C_q=q^2I^2_1+I^2_2+q^2I^2_3+q^{1/2}(1-q^2)I_1I_2I_3.
$$
As in the case of quantum algebras (see, for example, Chapter 6 in [2]),
at $q$ a root of unity this algebra has additional central elements.
\medskip

\noindent
{\bf Theorem 3.} {\it Let
$q^n=1$ for $n\in \N$ and $q^j \not=1$ for $0<j< n$. Then the elements
$$ C^{(n)}(I_j)=
\sum_{j=0}^{[\frac{n-1}{2}]} \kombd{n-j}{j} \frac{1}{n-j}
\Bigl( \frac{i}{q-q^{-1}} \Bigr)^{2j} I_j^{n-2j},\ \ \ j=1,2,3,
\eqno(16)
$$
belong to the center of $U'_q({\rm so}_3)$, where $[x]$ for $x
\in \R$ denotes the integral part of $x$.}
\medskip

The proof of this theorem is rather complicated (see [7]).
First it is proved that $C^{(n)}(I_1)$ belongs to the center of
$U'_q({\rm so}_3)$. This proof is
based on the formula $I_3 I_1^m=p_m(I_1)I_2+q_m(I_1)I_3$, where
$$\textstyle
p_m(x)=q^{-\frac{1}{2}} \bigl( \frac{x(q+q^{-1})}{2}
\bigr)^{m-1} \sum
\limits_{t=0}^{[\frac{m-1}{2}]} \komb{m}{2t+1} \bigl(
\bigl(\frac{q-q^{-1}}{q+q^{-1}}\bigr)^2-
\bigl(\frac{2}{x(q+q^{-1})}\bigr)^2 \bigr)^t,$$
$$\textstyle
q_m(x)=-q^{\frac{1}{2}} \frac{x(q-q^{-1})}{2} p_m(x)+
\bigl(\frac{x(q+q^{-1})}{2}\bigr)^m
\sum\limits_{t=0}^{[\frac{m}{2}]} \komb{m}{2t} \bigl(
\bigl(\frac{q-q^{-1}}{q+q^{-1}}\bigr)^2-
\bigl(\frac{2}{x(q+q^{-1})}\bigr)^2
\bigr)^t.$$
The proof also needs deep combinatorial identities, such that
$$\textstyle
\sum\limits_{t=0}^{[\frac{N-1}{2}]} \komb{N}{2t+1}
\komb{[\frac{N-1}{2}]-t}{[\frac{N-1}{2}]-C} \komb{t}{M}
=$$$$\textstyle=
4^{C-M} \komb{C}{M} (N-2C(1-N'))
\frac{(2[\frac{N}{2}])!
([\frac{N}{2}]+C-M)!
([\frac{N}{2}]-M)!}
{[\frac{N}{2}]!
(2C+1)!
(2[\frac{N}{2}]-2M)!
([\frac{N}{2}]-C)!},$$
$$\textstyle
\sum
\limits_{j=0}^{d} \komb{n-j}{j} \frac{2d-2j+1+(n-2d-1)^{n'}}{n-j}
\komb{[\frac{n-1}{2}]-d}{c-j}
\krat
\frac{(2[\frac{n}{2}]-2j)!
(n-1-c-d)!
([\frac{n}{2}]-c)!}
{([\frac{n}{2}]-j)!
(d-j+1-n')!
(2[\frac{n}{2}]-2c)!
(n-2d+n'-1)!}
=$$$$=
\frac{\komb{n-1-d}{d}\komb{n-1-c}{c}}
{\bigl(\frac{n}{2}\bigr)^{1-n'} (n-2d)^{n'}},
$$
where $N,n\in \N$, $0 \leq C,M\leq [\frac{N-1}{2}]$,
$0 \leq c,d \leq [\frac{n-1}{2}]$ and $N',n'=0$ or 1 such that
$N'= N \mod 2$ and $n'= n \mod 2$.
One also needs an extensive use of the fact that $q$ is a root of
unity.

If it is proved that $C^{(n)}(I_1)$ belongs to the center of
$U'_q({\rm so}_3)$, then we have to use the
automorphism
$\rho\!:\ U'_q({\rm so}_3) \rightarrow U'_q({\rm so}_3)$
defined by relations
$\rho(I_1)=I_2$, $\rho(I_2)=I_3$, $\rho(I_3)=I_1$.
This automorphism shows that $C^{(n)}(I_2)$ and $C^{(n)}(I_3)$
also belong to the center of $U'_q({\rm so}_3)$.
\medskip

\noindent
{\bf Conjecture 1.} {\it If $q$ is a root of unity as above, then
the elements $C_q$ and $C^{(n)}(I_j)$, $j=1,2,3$,
generate the center of $U'_q({\rm so}_3)$.}
\medskip

{\bf 4.} Central elements of the algebra $U'_q({\rm so}_m)$ for
any value of $q$ are found in [8] and [9]. They are given in the form of
homogeneous polynomials of elements of $U'_q({\rm so}_m)$.
If $q$ is a root of unity, then (as in the case of quantum algebras)
there are additional central elements of $U'_q({\rm so}_m)$.
\medskip

\noindent
{\bf Theorem 4.} {\it Let $q^n=1$ for $n \in \N$ and $q^j \not=1$
for $0<j< n$. Then the elements
$$
C^{(n)}(I^+_{kl})=
\sum_{j=0}^{[\frac{n-1}{2}]} \kombd{n-j}{j} \frac{1}{n-j}
\Bigl( \frac{i}{q-q^{-1}} \Bigr)^{2j} {I^+_{kl}}^{n-2j},\ \ \ k>l,
\eqno(17)
$$
belong to the center of $U'_q({\rm so}_m)$.}
\medskip

Let us prove this theorem for the algebra $U'_q({\rm so}_4)$ (for the
general case a proof is the same).
This algebra is generated by the elements
$I_{43}$, $I_{32}$, $I_{21}$. We introduce the elements
$I_{31}\equiv I^+_{31}$, $I_{42}\equiv I^+_{42}$,
$I_{41}\equiv I^+_{41}$ defined as indicated above.
Then the elements $I_{ij}$, $i>j$, satisfy the relations
$$
[I_{43},I_{21}]=0, \ \ \
[I_{32},I_{31}]_q=I_{21},\ \ \
[I_{21},I_{32}]_q=I_{31}, \eqno(18)
$$
$$
[I_{31},I_{21}]_q=I_{32},\ \ \
[I_{43},I_{42}]_q=I_{32},\ \ \
[I_{32},I_{43}]_q=I_{42},
$$
$$
[I_{42},I_{32}]_q=I_{43},\ \ \
[I_{31},I_{43}]_q=I_{41},\ \ \
[I_{21},I_{42}]_q=I_{41},
$$
$$
[I_{41},I_{21}]_q=I_{42},\ \ \
[I_{41},I_{31}]_q=I_{43}, \ \ \
[I_{42},I_{41}]_q=I_{21},
$$
$$
[I_{41},I_{32}]=0,\ \
[I_{43},I_{41}]_q=I_{31},\ \
[I_{42},I_{31}]=(q-q^{-1})(I_{21}I_{43}-I_{41}I_{32}) .
$$
If one wants to prove that an element
$X$ belongs to the center of
$U'_q({\rm so}_4)$, it is sufficient to prove that
$[X,I_{21}]=[X,I_{32}]=[X,I_{43}]=0$.

Let us consider the element $C^{(n)}(I_{21})$. It
belongs to the subalgebra $U'_q({\rm so}_3)$ generated by $I_{21}$,
$I_{32}$ and $I_{31}$:

\medskip
\qquad $I_{21}$ $I_{31}$ $I_{41}$

\qquad \phantom{$I_{22}$} $I_{32}$
\begin{picture}(0,0)(0,0)
\put(-2,-4){\line(0,1){27}}
\put(-2,23){\line(-1,0){34}}
\put(-36,23){\line(0,-1){15}}
\put(-36,8){\line(1,0){16}}
\put(-20,8){\line(0,-1){12}}
\put(-20,-4){\line(1,0){18}}
\end{picture}%
$I_{42}$

\qquad \phantom{$I_{23}$ $I_{33}$} $I_{43}$
\medskip

\noindent
It follows from Theorem 3 that $C^{(n)}(I_{21})$ commutes with
element $I_{32}$. Using the first relation in (18) we easily see
that $C^{(n)}(I_{21})$ commutes with
$I_{43}$ and therefore $C^{(n)}(I_{21})$ belongs to the center
of $U'_q({\rm so}_4)$.

Let us consider the element $C^{(n)}(I_{32})$. In
$U'_q({\rm so}_4)$ we separate two subalgebras
$U'_q({\rm so}_3)$:

\medskip
\qquad $I_{21}$ $I_{31}$ $I_{41}$

\qquad \phantom{$I_{22}$} $I_{32}$
\begin{picture}(0,0)(0,0)
\put(-2,-4){\line(0,1){27}}
\put(-2,23){\line(-1,0){34}}
\put(-36,23){\line(0,-1){15}}
\put(-36,8){\line(1,0){16}}
\put(-20,8){\line(0,-1){12}}
\put(-20,-4){\line(1,0){18}}
\end{picture}%
$I_{42}$

\qquad \phantom{$I_{23}$ $I_{33}$} $I_{43}$
\begin{picture}(0,0)(-2,2)
\put(-2,-4){\line(0,1){27}}
\put(-2,23){\line(-1,0){34}}
\put(-36,23){\line(0,-1){15}}
\put(-36,8){\line(1,0){16}}
\put(-20,8){\line(0,-1){12}}
\put(-20,-4){\line(1,0){18}}
\end{picture}%
\medskip

From Theorem 3 we have $[C^{(n)}(I_{32}),I_{21}]=
[C^{(n)}(I_{32}),I_{43}]=0$ and $C^{(n)}(I_{32})$ belongs to the center
of $U'_q({\rm so}_4)$.
A proof that the element $C^{(n)}(I_{43})$ belongs to the center
is the same as for $C^{(n)}(I_{21})$.

The elements $C^{(n)}(I_{31})$, $C^{(n)}(I_{42})$ and
$C^{(n)}(I_{41})$ belong to the center of $U'_q({\rm so}_4)$
since they belong to the sabalgebras
$U'_q({\rm so}_3)$ generated by triplets
$$ I_{41},\ \ I_{31},\ \ I_{43} \ \ \ \ {\rm and}\ \ \ \
I_{21},\ \ I_{41},\ \ I_{42}.
$$
(Note that $C^{(n)}(I_{31})$ and $C^{(n)}(I_{42})$ commute with $I_{42}$
and $I_{31}$, respectively, since $I_{42}=[I_{32},I_{43}]_q$ and
$I_{31}=[I_{21},I_{32}]_q$.) Theorem is proved.
\medskip

\noindent
{\bf Conjecture 2.} {\it If $q$ is a root of unity as above, then the
central elements of [9] and of Theorem 4 generate the center of
$U'_q({\rm so}_m)$.}
\medskip

{\bf 5.}
Let us consider the associative algebra
$U'_{q,\varepsilon}({\rm so}_3)$ (where $\varepsilon \ge 0$)
generated by three generators $J_1$, $J_2$, $J_3$ satisfying
the relations:
$$
[J_1,J_2]_q:=q^{1/2} J_1 J_2-q^{-1/2} J_2 J_1=J_3,\ \ \
[J_2, J_3]_q=J_1,\ \ \
[J_3, J_1]_q=\varepsilon^2 J_2.
$$
It is easily proved that
this algebra is isomorphic to the algebra
$U'_{q}({\rm so}_3)$ and the corresponding isomorphism is uniquely
defined by $J_1 \rightarrow
\varepsilon I_1$, $J_3 \rightarrow \varepsilon I_3$, $J_2
\rightarrow I_2$. Therefore, the elements
$$
\widetilde C^{(n)}(J_i,\varepsilon):=n\varepsilon^n
C^{(n)}(J_i/{\varepsilon}),\ i=1,3,\ \ \ \
\widetilde C^{(n)}(J_2,\varepsilon):=C^{(n)}(J_2)
$$
belong to the center of $U'_{q,\varepsilon}({\rm so}_3)$ if
$q^n=1$. By means of the contraction $\varepsilon \to 0$ we transform
the algebra $U'_{q,\varepsilon}({\rm so}_3)$ into the algebra
$U'_q({\rm iso}_2)$. Under this contraction the
commutation relations $[\widetilde
C^{(n)}(J_i,\varepsilon),J_k]=0$ transform into the
relations $[\widetilde C^{(n)}(J_i,0),J_k]=0$.
In other words, we have proved the following
\medskip

\noindent
{\bf Theorem 5.} {\it Let $q^n=1$ for $n \in \N$ and $q^j \not=1$ for
$0<j< n$. Then the elements $T_1^n$, $T_2^n$ and $C^{(n)}(I)$
belong to the center of the algebra $U'_q({\rm iso}_2)$.}
\medskip

It was shown in [10] that the element
$$
C_q=q^{-1}T^2_1+qT^2_2+q^{-3/2}(1-q^2)T_1T_2I
$$
is central in $U'_q({\rm iso}_2)$.
\medskip

\noindent
{\bf Conjecture 3.} {\it If $q$ is a root of unity as above, then
the elements $C_q$, $T_1^n$, $T_2^n$ and $C^{(n)}(I)$
generate the center of $U'_q({\rm iso}_2)$.}
\medskip

Using Theorem 5 and repeating the proof of Theorem 4 we prove the
following theorem:
\medskip

\noindent
{\bf Theorem 6.} {\it Let $q^n=1$ for $n \in \N$ and $q^j \not=1$ for
$0<j< n$. Then the elements
$$
C^{(n)}(I_{ij}),\ \ i>j,\ \ \ \ \ T_j^n,\ \ j=1,2,\cdots ,m,
$$
belong to the center of the algebra $U'_q({\rm iso}_m)$.}

\end{document}